\documentclass[a4paper,12pt]{amsart}
\usepackage{amssymb}
\usepackage[arrow,tips,matrix]{xy}
\usepackage{hyperref}

\calclayout

\theoremstyle{plain}

\newtheorem{theo}{Theorem}[section]
\newtheorem{prop}[theo]{Proposition}
\newtheorem{lemm}[theo]{Lemma}
\newtheorem{coro}[theo]{Corollary}

\theoremstyle{remark}
\newtheorem*{rema}{Remark}

\newcommand{\C}{\mathbb{C}}

\newcommand{\side}[1]{{}^{#1}\hskip-.075em}

\DeclareMathOperator{\GL}{GL} \DeclareMathOperator{\PGL}{PGL}
\DeclareMathOperator{\Tor}{Tor} 
\DeclareMathOperator{\Irr}{Irr} \DeclareMathOperator{\Hom}{Hom}

\begin{document}

\title{On extensions of representations for compact Lie groups}

\date{July 2, 2002. First draft on December 22, 1999.}

\author{Jin-Hwan Cho}
\address{School of Mathematics, Korea Institute for Advanced Study}
\email{chofchof@kias.re.kr}
\thanks{Jin-Hwan Cho would like to thank Osaka City University for its
hospitality during his visit when the first draft of the paper was
written.}

\author{Min Kyu Kim}
\address{Department of Mathematics, Korea Advanced Institute of Science and Technology}
\email{minkyu@math.kaist.ac.kr}

\author{Dong Youp Suh}
\address{Department of Mathematics, Korea Advanced Institute of Science and Technology}
\email{dysuh@math.kaist.ac.kr}
\thanks{Dong Youp Suh wishes to acknowledge the financial support
of the Korea Research Foundation made in the program year of 2001,
and Grant No.~R01-1999-00002 from the Interdisciplinary Research
Program of KOSEF}



\begin{abstract}
Let $H$ be a closed normal subgroup of a compact Lie group $G$
such that $G/H$ is connected. This paper provides a necessary and
sufficient condition for every complex representation of $H$ to be
extendible to $G$, and also for every complex $G$-vector bundle
over the homogeneous space $G/H$ to be trivial. In particular, we
show that the condition holds when the fundamental group of $G/H$
is torsion free.
\end{abstract}

\maketitle

\section{Introduction}

One of the classical problems in finite group theory is to
characterize extensions of representations. We mean an extension
of a representation in the following way: Given a normal subgroup
$H$ of a group $G$, a (complex) representation $\rho\colon
H\to\GL(n,\C)$ is called \emph{extendible to $G$} if there exists
a representation $\widetilde\rho\colon G\to\GL(n,\C)$ (called a
\emph{$G$-extension}) such that $\rho=\widetilde\rho$ on $H$. It
is to be noted that the dimension $n$ is not changed, since $\rho$
as a sub-representation is always contained in the restriction of
the induced representation of $\rho$ to $H$.

In the case of finite $G$, it is well known that every complex
irreducible representation of $H$, which is $G$-invariant under
conjugation (see Section~2 for the definition), is extendible to
$G$ if the second group cohomology $H^2(G/H,\C^*)$
vanishes~\cite[Theorem~11.7]{Isa76}. On the other hand the
extension problem for infinite groups has not been extensively
studied. In this article we study the problem for compact Lie
groups when $G/H$ is connected. Our main result is a necessary and
sufficient condition for every complex representation of $H$ to be
extendible to $G$. It is also shown that the condition is related
to a topological invariant, the fundamental group of $G/H$.

For any group $G$, let $G'$ denote the commutator subgroup of $G$.

\begin{theo} \label{theo:main_theorem}
Let $G$ be a compact Lie group and $H$ a closed normal subgroup
such that $G/H$ is connected. Then every complex representation of
$H$ is extendible to $G$ if and only if $H$ is a direct summand of
$G'H$.
\end{theo}

\begin{coro} \label{coro:main_corollary}
Let $G$ be a compact Lie group and $H$ a closed normal subgroup
such that $G/H$ is connected. Then every complex representation of
$H$ is extendible to $G$ if the fundamental group $\pi_1(G/H)$ is
torsion free, or equivalently if $(G/H)'$ is simply connected.
\end{coro}

Our theorem provides a complete characterization of the triviality
of complex $G$-vector bundles over the homogeneous space $G/H$.
Let $E$ be a complex $G$-vector bundle over $G/H$. We recall that
$E$ is \emph{trivial} if it is isomorphic to the product bundle
$G/H\times V$ for some complex $G$-module $V$. Since $E$ is
uniquely determined by the fiber at the identity element of $G/H$
(say $E_0$), the bundle $E$ is trivial if and only if $E_0$ as a
complex representation of $H$ is extendible to $G$.
Theorem~\ref{theo:main_theorem} leads us to the following
corollary.

\begin{coro}
Let $G$ be a compact Lie group and $H$ a closed normal subgroup
such that $G/H$ is connected. Then every complex $G$-vector bundle
over the homogeneous space $G/H$ is trivial if and only if $H$ is
a direct summand of $G'H$. \qed
\end{coro}

The existence of $G$-extensions plays an important role even in
equivariant $K$-theory. Let $X$ be a connected topological space
with a compact Lie group $G$ action. Let $H$ be the normal
subgroup of $G$ which consists of all elements of $G$ acting
trivially on $X$. Then the projection $G\to G/H$ induces the
canonical homomorphism $\phi\colon K_{G/H}(X)\to K_G(X)$ which
sends a $G/H$-vector bundle over $X$ to the same bundle viewed as
a $G$-vector bundle with the trivial $H$-action.

On the other hand, suppose that every complex irreducible
representation of $H$ is extendible to $G$. Then there is an
injective group homomorphism $e\colon R(H)\to R(G)$ between two
representation rings defined as follows.
For each irreducible complex $H$-module $U$ choose a $G$-extension
$U_G$, and define $e([U])=[U_G]$ where $[\ \ ]$ denote  the classes in the representation rings.
Then extend the definition of $e$ to $R(H)$ so that it defines a homomorphism $R(H)\to R(G)$.
For each complex $G$-module $V$ we can
associate the trivial complex $G$-vector bundle
$\underline{V}=X\times V$, which defines the natural homomorphism
$t\colon R(G)\to K_G(X)$.
We now define a group homomorphism
\begin{equation}\label{K-homomorphism}
\mu\colon R(H)\otimes K_{G/H}(X)\to K_G(X),\qquad (V,\xi)\mapsto t\circ e(V)\otimes \phi(\xi).
\end{equation}

This
homomorphism is an isomorphism. Indeed, the inverse is given as
follows. Let $\Irr(H)$ denote the set of all isomorphism classes
of complex irreducible representations of $H$. For each $[\chi]\in
\Irr(H)$ choose a $G$-extension of $\chi$, and let $V_\chi$ be the
corresponding $G$-module to the chosen $G$-extension. For a
complex $G$-vector bundle $E$ over $X$, the canonical isomorphism
\[
E\xrightarrow{\cong} \bigoplus_{[\chi]\in\Irr(H)}
\underline{V_\chi} \otimes \Hom_H(\underline{V_\chi},E)
\]
induces a group homomorphism $K_G(X)\to R(H)\otimes K_{G/H}(X)$
which is the desired inverse (see~\cite[Section~2]{CKMS99} for
more general arguments). Therefore we have a generalization of
Proposition~2.2 in~\cite{Seg68} which deals with the extreme case
when $G$ acts trivially on $X$.

%

\begin{coro}
Let $G$ be a compact Lie group and $H$ a closed normal subgroup
such that $G/H$ is connected. Let $X$ be a connected $G$-space
such that $H$ acts trivially on $X$. If $H$ is a direct summand
of $G'H$, then the map $\mu\colon R(H)\otimes K_{G/H}(X)\to K_G(X)$
in (\ref{K-homomorphism}) can be defined,  and it is a group isomorphism. \qed
\end{coro}


This article is organized as follows. In Section~2 we shall give
some basic notions and then show that a complex irreducible
representation of $H$, which is $G$-invariant under conjugation,
induces an associated projective representation of $G$ which may
be viewed as a $G$-extension in the projective representation
level. Section~3 is devoted to prove that every complex
representation of $H$ has a $G$-extension when $G/H$ is connected
and abelian. In Section~4 we shall proceed the study in the case
that $G/H$ is semisimple and connected. After showing that the
extension problem can be reduced to this case, we shall prove
Theorem~1.1.

The authors wish to thank Professor Mikiya Masuda of Osaka City
University for valuable discussions on the overall contents of the
article. The authors also wish to thank Professor I. Martin Isaacs
of University of Wisconsin and Professor Hi-joon Chae of Hong-Ik
University for helpful discussions on finite and Lie group
representations.

\section{Associated projective representations}

Let $G$ be a topological group and $H$ a closed normal subgroup of
$G$. By a (complex) \emph{representation} of $G$ we shall mean a
continuous homomorphism of $G$ into the general linear group
$\GL(n,\C)$ of nonsingular $n\times n$ matrices over the field
$\C$ of complex numbers. A representation $\rho\colon
H\to\GL(n,\C)$ is called \emph{extendible to $G$} if there exists
a representation $\widetilde\rho\colon G\to\GL(n,\C)$ (called a
\emph{$G$-extension} of $\rho$) such that
$\rho(h)=\widetilde\rho(h)$ for all $h\in H$.
\[
\SelectTips{cm}{}
\xymatrix{ H \ar[r]^-{\rho} \ar[d] & \GL(n,\C) \\
G \ar@{.>}[ur]^-{\widetilde\rho} }
\]
Moreover, it is enough to get a $G$-extension of $\rho$ that there
is a representation $\widetilde\rho\colon G\to\GL(n,\C)$ such that
its restriction to $H$ is isomorphic (or similar) to $\rho$, i.e.,
there exists a matrix $M\in\GL(n,\C)$ such that
$M^{-1}\widetilde\rho(h)M=\rho(h)$ for all $h\in H$.

Given a representation $\rho\colon H\to\GL(n,\C)$ the map
$\side{g}\chi\colon H\to\C$ defined by the conjugation
$\side{g}\chi(h)=\chi(g^{-1}hg)$ becomes a representation of $H$
for each $g\in G$. We say that $\rho$ is \emph{$G$-invariant} if
it is isomorphic to the conjugate representation $\side{g}\rho$
for all $g\in G$, which is a necessary condition of $\rho$ to be
extendible to $G$.

\medskip

In the following we assume that a representation $\rho\colon
H\to\GL(n,\C)$ is irreducible and $G$-invariant. Then there exists
a matrix $M_g\in\GL(n,\C)$ for each $g\in G$ such that
$M_g^{-1}\rho(h)M_g=\side{g}\rho(h)=\rho(g^{-1}hg)$ for all $h\in
H$. Since $\rho$ is irreducible, the Schur's lemma implies that
$M_g$ is unique up to multiplication by nonzero constant in
$\C^*=\C\setminus\{0\}$. So we are able to define a function
$\rho^*$ of $G$ into the projective linear group
$\PGL(n,\C)=\GL(n,\C)/\C^*$ by $\rho^*(g)=[M_g]$ for each $g\in
G$, where $[M_g]$ denotes the image of $M_g$ by the canonical
projection $\pi\colon\GL(n,\C)\to\PGL(n,\C)$.
\[
\SelectTips{cm}{}
\xymatrix{ H \ar[r]^-{\rho} \ar[d] & \GL(n,\C) \ar[d]^-{\pi} \\
G \ar[r]^-{\rho^*} & \PGL(n,\C) }
\]

\begin{lemm} \label{lemm:complex_projective_extension}
Let $G$ be a topological group and $H$ a compact normal subgroup
of $G$. Given a complex irreducible representation $\rho\colon
H\to\GL(n,\C)$ which is $G$-invariant, the function $\rho^*\colon
G\to\PGL(n,\C)$ defined above is a continuous homomorphism, called
the projective representation of $G$ associated with $\rho$.
Moreover, the image of $\rho^*$ is contained in
$U(n)/S^1\subset\PGL(n,\C)$ if $\rho$ is a unitary representation
of $H$.
\end{lemm}

\begin{proof}
It is immediate that $\rho^*$ is a homomorphism. Since $H$ is
compact we may assume that $\rho$ is a unitary representation of
$H$, i.e., the image of $\rho$ is contained in the unitary group
$U(n)$. Then $M_g$ is a constant multiple of a matrix in $U(n)$ so
that $\rho^*(g)$ is contained in $U(n)/S^1$ for all $g\in G$. For
the continuity of $\rho^*$ it suffices to show that the graph of
$\rho^*$ in $G\times\PGL(n,\C)$ is closed, since $U(n)/S^1$ is a
compact Hausdorff space.

Consider the family of continuous maps $\Phi_h\colon
G\times\GL(n,\C)\to\GL(n,\C)$ for each $h\in H$ given by
$(g,M)\mapsto\rho(h)M\rho(g^{-1}hg)^{-1}M^{-1}$. Then the set
\[
\bigcap_{h\in H}\Phi_h^{-1}(I)=\bigcup_{g\in G}\bigl\{(g,M)\in
G\times\GL(n,\C)\mid M\in\pi^{-1}(\rho^*(g)) \bigr\},
\]
is the inverse image of the graph of $\rho^*$ in
$G\times\PGL(n,\C)$ by the canonical projection $1\times\pi\colon
G\times\GL(n,\C)\to G\times\PGL(n,\C)$, which is obviously closed
in $G\times\GL(n,\C)$. Therefore the graph of $\rho^*$ is also
closed in $G\times\PGL(n,\C)$.
\end{proof}

We may say that $\rho$ is extendible to $G$ in the projective
representation level, since $\rho^*(h)=[\rho(h)]$ for all $h\in
H$, i.e., $\rho^*=\pi\circ\rho$ on $H$.
\[
\SelectTips{cm}{}
\xymatrix{ H \ar[r]^-{\rho} \ar[d] & \GL(n,\C) \ar[d]^-{\pi} \\
G \ar[r]_-{\rho^*} \ar@{.>}[ur]^-{\widetilde\rho} & \PGL(n,\C) }
\]
Note that any $G$-extension (if exists) $\widetilde\rho$ of $\rho$
is a lifting homomorphism of $\rho^*$, i.e.,
$\rho^*=\pi\circ\widetilde\rho$, since
$\rho^*(g)=[\widetilde\rho(g)]$ for all $g\in G$.

%
%
%

\begin{rema}
In case that $G$ is finite, choose a transversal $T$ containing
$e$ for $H$ in $G$ and set $M_e=I$, the identity matrix in
$\GL(n,\C)$. For each $t\in T$ and $h\in H$, the map $\rho'\colon
G\to\GL(n,\C)$ sending $th\mapsto M_t\rho(h)$ is a lifting (not
necessarily homomorphism) of $\rho^*$, i.e.,
$\pi\circ\rho'=\rho^*$, and it determines a cocycle $\beta$ in the
second group cohomology $H^2(G/H,\C^*)$, which depends only on
$\rho$. Moreover, $\rho$ is extendible to $G$ if and only if
$\beta$ is trivial, see~\cite[Theorem~11.7]{Isa76} for more
details.
\end{rema}

\section{Extensions when $G/H$ is connected abelian}

In this section we shall prove that every complex representation
of $H$ is extendible to $G$ when $G/H$ is compact, connected, and
abelian, that is a torus. We begin with a general result on
extensions of representations in the special case when $G=SH$ for
some closed subgroup $S$ of $G$.

\begin{lemm} \label{lemm:complex_extension_property}
Let $G$ be a compact topological group such that $G=SH$ for a
closed subgroup $S$ and a closed normal subgroup $H$ of $G$. Then
a complex representation $\rho\colon H\to\GL(n,\C)$ is extendible
to $G$ if and only if there exists a representation $\varphi\colon
S\to\GL(n,\C)$ such that
\begin{enumerate}
\item $\varphi=\rho$ on $S\cap H$, and
\item $\varphi(s)^{-1}\rho(h)\varphi(s)=\rho(s^{-1}hs)$ for all
$s\in S$ and $h\in H$.
\end{enumerate}
\end{lemm}

\begin{proof}
The necessity is obvious so we prove the sufficiency. Define a
function $\widetilde\rho\colon G\to\GL(n,\C)$ by
$\widetilde\rho(sh)=\varphi(s)\rho(h)$ for $s\in S$ and $h\in H$.
It is immediate that $\widetilde\rho=\rho$ on $H$. In this proof
we shall use the symbols $s,s'$ and $h,h'$ for elements in $S$
and $H$, respectively.

\emph{Claim: $\widetilde\rho$ is well-defined.\ } If $sh=s'h'\in
G$, then $(s')^{-1}s=h'h^{-1}\in S\cap H$. Then the condition (1)
implies that $\varphi(s')^{-1}\varphi(s)=\rho(h')\rho(h)^{-1}$
and thus
$\widetilde\rho(sh)=\varphi(s)\rho(h)=\varphi(s')\rho(h')=\widetilde\rho(s'h')$.

\emph{Claim: $\widetilde\rho$ is a homomorphism.\ } For
$sh,s'h'\in G$, the condition (2) implies that
\begin{align*}
\widetilde\rho((s'h')(sh))&=\varphi(s')\varphi(s)\rho(s^{-1}h's)\rho(h) \\
&=\varphi(s')\varphi(s)\varphi(s)^{-1}\rho(h')\varphi(s)\rho(h) \\
&=\widetilde\rho(s'h')\widetilde\rho(sh),
\end{align*}
since $(s'h')(sh)=(s's)(s^{-1}h's)h$ and $s^{-1}h's\in H$.

\emph{Claim: $\widetilde\rho$ is continuous.\ } The map $p\colon
S\times H\to G$ sending $(s,t)\mapsto st$ is a continuous
surjection. Since both $S$ and $H$ are compact, $p$ is a closed
map so that $G$ has the quotient topology induced by $p$.
\[
\SelectTips{cm}{}
\xymatrix{ S\times H \ar[d]_-{p} \ar[dr]^-{\widetilde\rho\circ p} \\
G \ar[r]^-{\widetilde\rho} & \GL(n,\C) }
\]
Then the continuity of $\widetilde\rho$ follows from the
universal property of the identification map $p$ since the
composition $\widetilde\rho\circ p\colon S\times H\to\GL(n,\C)$
sending $(s,t)\mapsto\varphi(s)\rho(t)$ is continuous.
\end{proof}

\begin{rema}
In case that $\rho$ is irreducible, the condition (2) in
Lemma~\ref{lemm:complex_extension_property} implies that $\varphi$
is a lifting homomorphism of the associated projective
representation $\rho^*$ (defined in the previous section) over
$S$, i.e., $\pi\circ\varphi=\rho^*$ on $S$. On the other hand, any
lifting homomorphism $\varphi$ of $\rho^*$ over $S$ satisfies the
condition (2).
\end{rema}

\medskip

Our main concern in this paper is to study extensions of
representations when $G$ is a compact Lie group and $H$ is a
closed normal subgroup of $G$ such that $G/H$ is connected. In
this case every complex representation $\rho$ of $H$ is
$G$-invariant. Indeed, for each $g\in G$, there is a continuous
path $g_t$ in $G$ from $g$ to an element $h\in H$ since every
connected component of $G$ contains an element of $H$. Then the
path $g_t$ induces a continuous family of conjugate
representations $\side{g_t}\rho$ so that all representations
$\side{g_t}\rho$ are isomorphic (see~\cite[Lemma~38.1]{CoFl64} for
more general result). In particular, $\side{g}\rho=\side{g_0}\rho$
and $\rho=\side{h}\rho=\side{g_1}\rho$ are isomorphic.


Let $\rho$ be a complex irreducible representation of $H$. Since
$\rho$ is always $G$-invariant, the associated projective
representation $\rho^*$ exists by
Lemma~\ref{lemm:complex_projective_extension}. To get a
$G$-extension of $\rho$ we shall first find a closed subgroup $S$
of $G$ such that $G=SH$, and then construct a lifting homomorphism
$\varphi$ of $\rho^*$ over $S$ (so that the condition (2) is
satisfied). Finally modifying $\varphi$ a little to satisfy the
condition (1) we may get a $G$-extension of $\rho$.

\begin{lemm} \label{lemm:splitting_circle}
Let $G$ be a compact Lie group and $H$ a closed normal subgroup
such that $G/H\cong S^1$. Then there exists a circle subgroup $S$
of $G$ such that $G=SH$ and $S\cap H$ is finite cyclic.
\end{lemm}

\begin{proof}
Let $G_0$ denote the identity component of $G$. Since the
canonical projection $p\colon G\to G/H$ is open and closed,
$p(G_0)$ is a connected component of $G/H$ so that $p(G_0)=G/H$.
It is well known in Lie group theory~\cite[Theorem~6.15]{HoMo98}
that $G_0=Z_0G_0^{\prime}$, where $Z_0$ is the identity component
of the center of $G_0$, which is a torus and $G_0^{\prime}$ is the
commutator subgroup of $G_0$. Then $G_0^{\prime}\subset G_0\cap
H\subset H$ since $G/H=G_0/(G_0\cap H)$ is abelian, and thus
$p(Z_0)=G/H$. Using the isomorphism $G/H\cong U(1)$ we may view
$p|_{Z_0}$ as a one-dimensional unitary representation of the
torus $Z_0$. It is elementary in representation theory that there
exists a circle subgroup $S\subset Z_0$ such that $p(S)=G/H$.
Therefore $G=SH$ and, furthermore, the proper subgroup $S\cap H$
of the circle group $S$ is finite cyclic.
\end{proof}

\begin{lemm} \label{lemm:splitting_circle_two}
Let $T$ be a maximal torus in $U(n)$. Then the exact sequence
$0\to S^1\to T\to T/S^1\to 0$ splits. Here $S^1$ is identified
with the subgroup of $U(n)$ consisting of constant multiples $zI$
for $z\in S^1\subset \C$ where $I$ denotes the identity matrix.
\end{lemm}

\begin{proof}
Since any maximal torus $T$ in $U(n)$ is conjugate to the subgroup
$\Delta(n)\subset U(n)$ of diagonal matrices
\[
D(z_1,\dotsc,z_n)=\begin{pmatrix} z_1 & & \\ & \ddots & \\ & & z_n
\end{pmatrix}, \quad z_i \in S^1,
\]
it suffices to show that the exact sequence $0\to S^1\to
\Delta(n)\to \Delta(n)/S^1\to 0$ splits.
But the
splitting is immediate because of the homomorphism $\Delta(n)\to
S^1$ mapping a diagonal matrix $D(z_1,\dotsc,z_n)$ to the constant
multiple $z_1I\in S^1$.
\end{proof}

\begin{prop} \label{prop:complex_circle_extension}
Let $G$ be a compact Lie group and $H$ a closed normal subgroup
such that $G/H\cong S^1$. Then every complex representation of $H$
is extendible to $G$.
\end{prop}

\begin{proof}
Let $\rho\colon H\to\GL(n,\C)$ be a given representation. Since
$H$ is compact,
we may assume that all the images of $\rho$ are contained in
$U(n)\subset\GL(n,\C)$. Moreover, it is enough to prove the case
that $\rho$ is irreducible. Since $G/H\cong S^1$ is connected,
$\rho$ is $G$-invariant so that the associated projective
representation $\rho^*\colon G\to
U(n)/S^1\subset\PGL(n,\C)=\GL(n,\C)/\C^*$ exists by
Lemma~\ref{lemm:complex_projective_extension}. From
Lemma~\ref{lemm:splitting_circle} we can choose a circle subgroup
$S$ of $G$ such that $G=SH$ and $S\cap H$ is finite cyclic.

We shall find a lifting homomorphism $\varphi_0\colon S\to U(n)$
of $\rho^*$ over $S$. Since $\rho^*(S)$ is compact, connected, and
abelian, it is a torus in $U(n)/S^1$. Note that every maximal
torus in $U(n)/S^1$ has the form $T/S^1$ for some maximal torus
$T$ of $U(n)$~\cite[Theorem~2.9, Chapter~IV]{BrtD85}. Choose a
maximal torus $T$ of $U(n)$ such that $\rho^*(S)\subset T/S^1$. By
Lemma~\ref{lemm:splitting_circle_two} the exact sequence $0\to
S^1\to T\xrightarrow{\pi} T/S^1\to 0$ splits, i.e., the canonical
projection $\pi\colon T\to T/S^1$ has a continuous section
(homomorphism) $s\colon T/S^1\to T$ such that the composition
$\pi\circ s$ is the identity map of $T/S^1$. Then
$\varphi_0=s\circ\rho^*|_S$ is a desired lifting homomorphism of
$\rho^*$ over $S$.
\[
\SelectTips{cm}{}
\xymatrix{ & **[r] \hphantom{/}T\hphantom{S^1}\subset U(n) \ar@<1ex>[d]^-{\pi} \\
S \ar[r]_-{\rho^*} \ar@{.>}[ur]^-{\varphi_0} & **[r] T/S^1\subset
U(n)/S^1 \ar[u]^-{s} }
\]

Let $t_0$ denote a generator of the finite cyclic group $S\cap H$.
Since $\pi\circ\varphi_0=\rho^*=\pi\circ\rho$ on $S\cap H$,
$\varphi_0(t_0)=\xi\rho(t_0)$ for some constant $\xi\in S^1\subset
\C^*$. Note that $\xi$ is an $n$-th root of unity, where $n$ is
the order of $S\cap H$. So it is possible to choose a
one-dimensional unitary representation $\tau$ of the circle group
$S$ such that $\tau(t_0)=\xi^{-1}$. Then the unitary
representation $\varphi=\tau\otimes\varphi_0$ satisfies the
conditions (1) and (2) in
Lemma~\ref{lemm:complex_extension_property}.
\end{proof}

\begin{coro} \label{coro:complex_abelian_extension}
Let $G$ be a compact Lie group and $H$ a closed normal subgroup
such that $G/H$ is connected and abelian. Then every complex
representation of $H$ is extendible to $G$.
\end{coro}

\begin{proof}
Since $G/H$ is compact, connected, and abelian, it is isomorphic
to a torus. So we have a finite chain of subgroups
\[
H=H_0 \vartriangleleft H_1 \vartriangleleft \dotsb
\vartriangleleft H_{n-1} \vartriangleleft H_n=G
\]
such that $H_i$ is normal in $H_{i+1}$ and $H_{i+1}/H_i\cong S^1$.
Applying Proposition~\ref{prop:complex_circle_extension}
inductively, any representation of $H$ is extendible to $G$.
\end{proof}

\section{Extensions when $G/H$ is connected}

In this section we consider the general case, so $G/H$ will be
assumed to be connected (not necessarily abelian).
In this case the commutator subgroup $(G/H)'=G'H/H$ of $G/H$ is
semisimple connected~\cite[Theorem~6.18]{HoMo98}. The following
proposition reduces the extension problem to the case that $G/H$
is semisimple and connected.

\begin{prop} \label{prop:complex_connected_extension_condition}
Let $G$ be a compact Lie group and $H$ a closed normal subgroup of
$G$ such that $G/H$ is connected. A complex representation of $H$
is extendible to $G$ if and only if it is extendible to $G'H$.
\end{prop}

\begin{proof}
The necessity is obvious, and the sufficiency follows from
Corollary~\ref{coro:complex_abelian_extension} since the factor
group $G/G'H\cong(G/H)/(G'H/H)=(G/H)/(G/H)'$ is compact,
connected, and abelian, that is a torus.
\end{proof}

In the case that $G/H$ is semisimple connected, the following
result is well known in Lie group theory (see for
instance,~\cite[Proposition~6.14]{HoMo98}).

\begin{lemm} \label{lemm:semisimple_structure}
Let $G$ be a compact Lie group and $H$ a closed normal subgroup
such that $G/H$ is semisimple and connected. Then there is a
semisimple connected closed normal subgroup $S$ in $G$ such that
$G=SH$ and the map $S\times H\to G$ sending $(s,h)\mapsto sh$ is a
homomorphism with a discrete kernel isomorphic to $S\cap H$. \qed
\end{lemm}

\begin{rema}
Proposition~6.14 in~\cite{HoMo98} deals with the case when $G$ is
connected. However, the same proof holds even if $G$ is not
connected, since $G/H$ is connected. Moreover, we can find the
fact in the proof that $S$ is semisimple and connected.
\end{rema}

The following result implies that the existence of a $G$-extension
when $G/H$ is semisimple and connected is completely determined by
the restriction of a given representation to $S\cap H$.

\begin{prop} \label{prop:complex_semisimple_extension}
Under the hypotheses of Lemma~\ref{lemm:semisimple_structure}, a
complex irreducible representation $\rho$ of $H$ is extendible to
$G$ if and only if $\rho$ is trivial on $S\cap H$, i.e.,
$\rho(g)=I$, the identity matrix, for all $g\in S\cap H$.
\end{prop}

\begin{proof}
It is immediate that $S$ commutes with $H$, since the map
$S\times H\to G$ sending $(s,h)\mapsto sh$ is a homomorphism. To
prove the sufficiency, it is enough to choose the trivial
representation $\varphi$ of $S$, i.e., $\varphi(s)=I$ for all
$s\in S$. Since $S$ commutes with $H$, the two conditions (1) and
(2) in Lemma~\ref{lemm:complex_extension_property} are satisfied
immediately.

On the other hand, suppose $\widetilde\rho$ is a $G$-extension of
$\rho$. Since $S$ commutes with $H$, we have
$\widetilde\rho(s)^{-1}\rho(h)\widetilde\rho(s)=\rho(h)$ for all
$s\in S$ and $h\in H$. Then the Schur's lemma implies that
$\widetilde\rho(s)$ is constant for all $s\in S$, so we may view
the restriction $\widetilde\rho|_S$ as a one-dimensional complex
representation of $S$. Since semisimple Lie groups have no
nontrivial abelian factor group, the trivial representation is
the unique one-dimensional complex representation of $S$.
Therefore, $\widetilde\rho$ is trivial on $S$, in particular, on
$S\cap H$.
\end{proof}

\begin{rema}
Note that the number of $G$-extensions (if exist) is exactly one,
since every $G$-extension should be trivial on $S$.
\end{rema}

\begin{coro} \label{coro:complex_extension_property_semisimple}
Let $G$ be a compact Lie group and $H$ a closed normal subgroup
such that $G/H$ is semisimple and connected. Every complex
representation of $H$ is extendible to $G$ if and only if $H$ is a
direct summand of $G$, i.e., $G\cong S\times H$ for some subgroup
$S$ of $G$.
\end{coro}

\begin{proof}
The sufficiency is obvious so we prove the necessity. If $H$ is
not a direct summand of $G$, then $S\cap H$ in
Lemma~\ref{lemm:semisimple_structure} contains a nontrivial
element, say $s_0$. Since a faithful representation of $H$ always
exists~\cite[Theorem~4.1, Chapter~III]{BrtD85}, we can choose an
irreducible sub-representation $\rho$ of $H$ such that $\rho(s_0)$
is not trivial. Then $\rho$ does not extend to a representation of
$G$ by Proposition~\ref{prop:complex_semisimple_extension}.
\end{proof}

We shall now prove the main result in this paper. For the second
statement of Theorem~\ref{theo:main_theorem}, we need the
following lemma giving a relation between the normal subgroup
$S\cap H$ in Lemma~\ref{lemm:semisimple_structure} and the
fundamental group of $G/H$.

\begin{lemm} \label{lemm:surjective_homomorphism}
Under the hypotheses of Lemma~\ref{lemm:semisimple_structure},
there exists a surjective homomorphism $\pi_1(G/H)\to S\cap H$.
\end{lemm}

\begin{proof}
Since $S/(S\cap H)=G/H$, the restriction of the canonical
projection $p\colon G\to G/H$ on $S$ is surjective and its kernel
$S\cap H$ is discrete. It follows that $p|_S$ is a covering
homomorphism of $G/H$. From the uniqueness of the universal
covering homomorphism $\widetilde q\colon\widetilde{G/H}\to G/H$,
there exists a covering homomorphism $q\colon\widetilde{G/H}\to S$
such that the diagram
\[
\xymatrix{\widetilde{G/H} \ar[d]_-{\widetilde q} \ar[r]^-{q} & S
\ar[dl]^-{p|_S} \\
G/H }
\]
commutes (compare with~\cite[Proposition~9.12]{HoMo98}). Since
$S\cap H=\ker p|_S=q(\ker\widetilde q)$ and $\ker\widetilde q$ is
isomorphic to $\pi_1(G/H)$, we have a surjective homomorphism of
$\pi_1(G/H)$ onto $S\cap H$.
\end{proof}

\renewcommand{\thetheo}{1.1~\textnormal{(rephrased)}}
\begin{theo}
Let $G$ be a compact Lie group and $H$ a closed normal subgroup
such that $G/H$ is connected. Then every complex representation of
$H$ is extendible to $G$ if and only if $H$ is a direct summand of
$G'H$.
\end{theo}

\begin{proof}
Since the factor group $G'H/H=(G/H)'$ is semisimple and connected,
the theorem follows immediately from
Proposition~\ref{prop:complex_connected_extension_condition} and
Corollary~\ref{coro:complex_extension_property_semisimple}.
\end{proof}

\begin{proof}[Proof of Corollary~\ref{coro:main_corollary}]
We claim that $\Tor(\pi_1(G/H))$, the torsion subgroup of
$\pi_1(G/H)$, is isomorphic to $\pi_1((G/H)')$. Denote by $T$ the
torus $(G/H)/(G/H)'$. Then the homotopy exact sequence of the
fibration $(G/H)'\to G/H\to T$ implies that
$\pi_1(G/H)\cong\pi_1((G/H)')\oplus\pi_1(T)$, since the second
homotopy group of a compact Lie group vanishes,
see~\cite[Proposition~7.5, Chapter~V]{BrtD85}. Since $(G/H)'$ is
semisimple, $\pi_1((G/H)')$ is finite~\cite[Remark~7.13,
Chapter~V]{BrtD85} so that it is isomorphic to $\Tor(\pi_1(G/H))$
as we claimed. Therefore, the condition of $\pi_1(G/H)$ being
torsion free is equivalent to $(G/H)'$ being simply connected.

By Lemma~\ref{lemm:semisimple_structure}
and~\ref{lemm:surjective_homomorphism}, $G'H=SH$ for some
semisimple connected closed normal subgroup $S$ in $G'H$ and there
is a surjective homomorphism $\pi_1(G'H/H)=\pi_1((G/H)')\to S\cap
H$. Therefore, if $(G/H)'$ is simply connected, then
$\pi_1((G/H)')=S\cap H=\{e\}$ so that $H$ is a direct summand of
$G'H$.
\end{proof}

\providecommand{\bysame}{\leavevmode\hbox
to3em{\hrulefill}\thinspace}

\end{document}